\documentclass[10pt]{amsart}

\usepackage{latexsym}
\usepackage{amssymb}
\usepackage{enumerate}
\usepackage[all]{xy}

\newif\ifpdf
\ifx\pdfoutput\undefined
	\pdffalse
\else
	\pdfoutput=1
	\pdftrue
\fi
\ifpdf
	\usepackage[pdftex]{graphicx}
	\DeclareGraphicsExtensions{.pdf}
	\usepackage{hyperref}
\else
	\usepackage{graphicx}
	\DeclareGraphicsExtensions{.eps}
	\usepackage{hyperref}
\fi

\DeclareSymbolFont{EulerScript}{U}{eus}{m}{n}
\DeclareSymbolFontAlphabet\mathscr{EulerScript}

\newcounter{Tlistc}

\newtheoremstyle{mythm}{10pt}{10pt}{\it}{}{\bf}{.}{ }{}
\theoremstyle{mythm}
\newtheorem{prop}{Proposition}[section]
\newtheorem{lemma}[prop]{Lemma}
\newtheorem{cor}[prop]{Corollary}
\newtheorem{thm}[prop]{Theorem}

\newtheoremstyle{myex}{10pt}{10pt}{\rm}{}{\bf}{.}{ }{}
\theoremstyle{myex}
\newtheorem*{ack}{Acknowledgment}

\newtheorem{rem}[prop]{Remark}
\newtheorem{definition}[prop]{Definition}

\newcommand{\Aa}{\mathbb A}
\newcommand{\Bb}{\mathbb B}
\newcommand{\Dd}{\mathbb D}
\newcommand{\Ff}{\mathbb F}
\newcommand{\Pp}{\mathbb P}
\newcommand{\Zz}{\mathbb Z}

\newcommand{\Cc}{\mathbb C}

\newcommand{\Qq}{\mathbb Q}
\newcommand{\cC}{\mathcal C}
\newcommand{\cD}{\mathcal D}
\newcommand{\eps}{\varepsilon}

\newcommand{\Ho}{H}
\newcommand{\Id}{\hbox{Id}}
\newcommand{\Sing}{\hbox{Sing}}

\newcommand{\GL}{\hbox{GL}}

\newcommand\enet[1]{\renewcommand\theenumi{#1}
\renewcommand\labelenumi{\theenumi}}

\title[Twisted Alexander Polynomials...]{Twisted Alexander polynomials\\
of Plane Algebraic Curves}

\author[J.I. Cogolludo]{J.I. Cogolludo}
\address{Departamento de Matem\'aticas\\
Universidad de Zaragoza}
\email{jicogo@unizar.es}

\author[V.Florens]{V. Florens}
\address{Departamento de \'Algebra y Geometr{\'\i}a\\
Universidad de Valladolid}
\email{vincent$\_$florens@yahoo.fr}

\thanks{First author is partially supported by MTM2004-08080-C02-02. Second author
is supported by MCHF-2001-0615.}

\begin{document}

\begin{abstract}
We consider the Alexander polynomial of a plane algebraic curve twisted by a linear 
representation. We show that it divides the product of the polynomials of the 
singularity links, for unitary representations. Moreover, their quotient is given 
by the determinant of its Blanchfield intersection form. Specializing in the 
classical case, this gives a geometrical interpretation of Libgober's divisibility 
Theorem. We calculate twisted polynomials for some algebraic curves and show how 
they can detect Zariski pairs of equivalent Alexander polynomials and that they are 
sensitive to nodal degenerations.
\end{abstract}

\maketitle

\section{Introduction}
Zariski~\cite{Za} used the fundamental group of the complement of a plane algebraic 
curve to show that there exist sextics with the same combinatorics (degree of 
irreducible components, local type of singularities,...) but different embeddings.
This comes to show that the position of singularities has an effect on the topology 
of the curve, and that the fundamental group is sensitive to these phenomena. Still,
since there is no classification of fundamental groups of curves, and the isomorphism
problem is undecidable, one cannot directly use this invariant in an effective way.
The Alexander polynomial is more manageable and also sensitive to the
position of singularities, see~\cite{Ar,AC,De,Li1,Oka1}.
Libgober~\cite{Li1} showed in particular that it divides the product of the local
polynomials, associated with its singular points. This result was sharpened by Degtyarev~\cite{De2} which described the type of singular points that may contribute 
to the global polynomial. Thereafter other invariants of the Alexander modules have
been considered in order to produce Zariski pairs with non isomorphic fundamental 
groups which cannot be distinguished by the Alexander polynomial. This is the case for
the Characteristic Varieties first considered by Libgober~\cite{Li2} and 
others~\cite{ACC,ACCT},
or the $\theta$-polynomials considered by Oka~\cite{Oka2}. Invariants related to
representations of the fundamental group such as the number of epimorphisms on finite groups
(also called Hall invariants)~\cite{AC}, or existence of dihedral coverings~\cite{ACCT}
have also been used in this context. Furthermore, the study of lower central 
series~\cite{ACCM,R} play an important role in the case of fundamental groups of
line arrangements. 

In knot theory, a strategy to study problems that the Alexander polynomial is not strong
enough to solve is to consider non-Abelian invariants, twisted by a linear representation 
of the fundamental group -see, for instance~\cite{Ki,Lin,Wa}. For mutation and 
concordance questions, Kirk and Livingston~\cite{KL1,KL2} recently developed their properties in the general context of CW-complexes. Their connection with lower central
series, though reasonable~\cite{Coc}, is far from being well understood.

In this paper, we use their results to establish the relationship between the Alexander polynomial of a plane algebraic curve, twisted by a unitary representation, and the product of the local ones.

For a finite CW-complex $X$ and a finite dimensional vector space $V$ over $\Ff$, one
defines the Alexander polynomials $\Delta_{\varepsilon,\rho}^i(X)$  in $\Ff[t^{\pm 1}]$, associated with a representation $\rho: \pi_1(X) \rightarrow \GL(V)$ and 
a choice of an epimorphism $\varepsilon : \pi_1(X) \rightarrow \Zz$.
Twisted Alexander polynomials of links are obtained by applying the invariants to their
exterior in $S^3$. In particular, Wada~\cite{Wa} considered the invariant
$\Delta_{\varepsilon,\rho}(L)= \Delta_{\varepsilon,\rho}^1(L) /
\Delta_{\varepsilon,\rho}^0(L)$.
Given an affine algebraic curve $\mathcal C$, with
exterior $X$ in a sufficiently large ball $\Bb^4$, one defines
$\Delta_{\varepsilon,\rho}^*(\mathcal C)= \Delta_{\varepsilon,\rho}^*(X)$ and sets
$\Delta_{\varepsilon,\rho}(\mathcal C)= 
\Delta_{\varepsilon,\rho}^1(\mathcal C)/ \Delta_{\varepsilon,\rho}^0(\mathcal C)$. 
The related Blanchfield~\cite{Bl} intersection form $\varphi^{\eps,\rho}(\cC)$ is 
defined on $H_2^{\varepsilon,\rho}(X;\Ff[t^{\pm 1}]$).

\vspace*{14pt}
In the formulas presented here, equalities should be understood up to units 
in~$\Ff[t^{\pm 1}]$.

\begin{thm}
\label{thm-main-intro}
Let $\cC=\cC_1\cup \dots \cup \cC_r$ be an affine algebraic curve with
$r$ irreducible components. Consider an epimorphism 
$\varepsilon: \pi_1(X) \rightarrow \Zz$ and a unitary representation
$\rho : \pi_1(X) \rightarrow \GL(V)$, such that $\Ff$ is a subfield of $\mathbb{C}$, 
closed by conjugation. 

Let $s=\# \Sing(\cC)$ and consider $(S^3_k,L_k)$ the local spheres and link 
singularities, $k=1,...,s$ plus the link at infinity $(S^3_\infty,L_\infty)$, where
$S^3_\infty=\partial \Bb^4$. 
Let us also denote by $(\varepsilon_k,\rho_k)$ 
the induced representations on $\pi_1(S^3_k\setminus L_k)$, 
$k=1,...,s,\infty$.
If the Alexander modules of $L_k$ are torsion for any $k$, then
$$\alpha \cdot \prod_{k=1,...,s,\infty} \Delta_{\varepsilon_k,\rho_k}(L_k) =
\overline{\Delta}_{\varepsilon,\rho}(\mathcal C) \cdot 
\Delta_{\varepsilon,\rho}(\mathcal C)
\cdot \det\varphi^{\varepsilon,\rho}(\cC),$$
where $\alpha=\prod_{\ell=1}^r 
\det(\Id - \rho(\nu_\ell) t^{\varepsilon(\nu_\ell)})^{s_\ell-\chi(\mathcal C_\ell)}$,
with $s_\ell=\# \Sing (\cC) \cap \cC_\ell$, and $\nu_\ell$ a meridian of $\cC_\ell$.
\end{thm}

Consider now the case of $(\varepsilon,\rho)$ such that $\rho(x)=1 \in \Qq$ for all 
$x \in \pi_1(X)$, and $\varepsilon$ sends all the meridians of $\mathcal C$ to 
$1 \in \Zz$ (or $t$ if denoted multiplicatively). We obtain the classical 
Alexander polynomial and Theorem~\ref{thm-main-intro} provides the following
geometrical interpretation of Libgober's Divisibility Theorem~\cite{Li1}.

\begin{cor}
Let $\Delta_{\mathcal C}$ be the classical Alexander polynomial of $X$, and 
$\Delta_{L_k}$ be the local Alexander polynomials. If $\varphi^t(\cC)$ is the intersection 
form with twisted coefficients in $\Qq[t^{\pm 1}]$, then
$$(t-1)^{1-\chi(\cC)}\prod_{k=1,...,s,\infty} \Delta_{L_k} = \Delta_{\mathcal C}^2 \cdot \det \varphi^t(\cC).$$
\end{cor}

Going back to knot theory, Milnor~\cite{Mi1,Mi2} showed that 
the Alexander polynomial essentially coincides with
the Franz-Reidemeister torsion of the link complement. Turaev~\cite{Tu1,Tu2}
further developed this construction, which provided new proofs for
several classical results.
The first version of twisted Alexander polynomial for knots
was due to Lin~\cite{Lin}. Wada~\cite{Wa} generalized it as an 
invariant of finitely presentable groups endowed with
a representation, in terms of Fox calculus. Then Kitano~\cite{Ki} showed
that it coincides with a torsion of the knot complement,
in the acyclic case.
Kirk and Livingston~\cite{KL1} extended his construction to the non-acyclic case,
and any finite CW-complex.

On the other hand, a procedure to compute the group of a plane curve was developed by
Zariski~\cite{Za} and Van Kampen~\cite{VK} and expressed by Moishezon~\cite{Mo} in terms 
of braid monodromy. Libgober~\cite{Li2} showed that the complement of the (affine) curve 
has the homotopy type of the $2$-dimensional complex corresponding to this presentation 
of the group. From this the relationship between the twisted Alexander polynomial 
and the related torsion can be established, similar to those in knot theory. The 
duality theorem, due to Kirk and Livingston~\cite{KL1} in this context, 
is applied in order to obtain the main result, see Theorem~\ref{main}.

We tried to write a paper as self-contained as possible. All the basics involved 
here and the duality theorem can be found in~\cite{KL1}. In Section~\S\ref{tw} we recall
basic constructions on the Reidemeister torsion of a CW-complex, in the non-acyclic case, and the relations with twisted Alexander polynomials. In Section~\S\ref{dua} we recall
the duality theorem for twisted torsion, due to~\cite{Mi2} in
the classical case and to~\cite{KL1} in the twisted case. In Section~\S\ref{hist},
we give a brief historical approach on the twisted torsion
and Alexander polynomials of links. Section~\S\ref{planecurves} is devoted
to plane curves, and to the proof of both Theorem~\ref{main}
and Corollary~\ref{corro}. Also we show the connection between Characteristic Varieties 
and twisted Alexander polynomials. Finally, Section~\S\ref{exs} illustrates
some examples of curves with trivial classical Alexander polynomials, but 
non-trivial twisted Alexander polynomials, using both abelian and non-abelian
representations.

Note finally that a several-variable version of twisted polynomials could have been considered, for an $\varepsilon$ related to the universal Abelian covering. For 
technical reasons, since $\Ff[t^{\pm 1}]$ is a principal ideal domain, we restrict 
ourselves to this case.

\begin{ack}
This paper was written while the second author visited the Universidad 
de Zaragoza. He appreciates its motivating mathematical atmosphere. We both thank
E.Artal for his support and useful suggestions.
\end{ack}

\section{Twisted Alexander polynomials}
\label{tw}

\subsection{Torsion of a chain complex}

Let $C_*$ be a finite chain complex:
$$ C_* = C_n \stackrel{\partial} \longrightarrow \dots  \stackrel{\partial} 
\longrightarrow C_0,$$
where $C_i$ are finite dimensional $\Ff$-vector spaces, and $\partial \circ \partial =0$.
Choose a basis $c_i$ for $C_i$, $\bar h_i$ a basis for the homology $\Ho_i(C_*)$ and 
$h_i$ a lift of $\bar h_i$ in $C_i$. Note that if $\Ho_i(C_*)=0$ for all $i$,
the complex is called \emph{acyclic}. Let $b_i$ be a basis of the image of 
$\partial : C_{i+1} \longrightarrow C_i$ and $\widetilde{b}_i$ be a lift of $b_i$ in 
$C_{i+1}$. One easily checks that $b_i h_i \widetilde{b}_{i-1}$ is a basis of $C_i$.
Denote by $[u,v]$ the determinant of the transition matrix from $u$ to $v$.

\begin{definition}
The  \emph{torsion} of $(C_*;c,h)$ is
$$\tau(C_*;c,h) =   \prod_{i=0}^n [b_{i} 
h_{i} \widetilde{b}_{i-1} | c_{i}]^{(-1)^{i+1}} \in \Ff^* / \{\pm 1\}.$$
\end{definition}

Note that in the literature the torsion
is sometimes defined as the inverse of $\tau(C_*;c,h)$.
The torsion does not depend on the choice of $b$ and its lifts.
It depends on the choice of $c$ and $h$ as follows:
$$ \tau(C_*; c^\prime,h^\prime) = \tau(C_*; c, h) \prod_i  
\Big( \frac{[h^\prime_{i} | h_i]}{[c^\prime_{i} | c_{i}]} \Big)^{(-1)^{i+1}}.$$

The following classical lemma is very useful for computations. 

\begin{lemma}[\cite{Mi2}]
\label{mult}
Consider a short exact sequence of complexes 
$$0 \longrightarrow C^\prime \longrightarrow C \longrightarrow 
C^{\prime \prime} \longrightarrow 0,$$
where the complexes and their homology are based, with compatible bases.
Let $\mathcal H$, with torsion $\tau(\mathcal H)$, be the related long exact 
sequence in homology, viewed as a based acyclic complex. One has
$$\tau(C) = \tau(C^\prime) \tau(C^{\prime \prime}) \tau(\mathcal H).$$
\end{lemma}

\subsection{Twisted chain complexes}

In this section, $X$ is a finite CW-complex, with $\pi=\pi_1(X,x)$ for $x \in X$. 
Let us fix an epimorphism 
$$\varepsilon : \pi  \longrightarrow \Zz.$$
Note that $\eps$ extends naturally to an epimorphism of algebras
$\varepsilon: \Zz[\pi] \rightarrow \Zz[\Zz]$, which will be also denoted by $\eps$.
We identify $\Zz[\Zz]$ with $\Zz[t^{\pm 1}]$.
Consider now an $\Ff$-vector space $V$ of finite dimension and a representation
$$\rho: \pi \longrightarrow \GL(V).$$ 
If $\widetilde{X} \rightarrow X$ denotes the universal covering, the cellular chain 
complex $C_*(\widetilde{X}; \Ff)$ is an $\Ff[\pi]$-module generated by the lifts of 
the cells of $X$. Consider the $\Ff[\pi]$-module $\Ff[t^{\pm 1}] \otimes_{\Ff} V$, 
where the action is induced by $\varepsilon \otimes \rho$, as follows:
$$ (p \otimes v) \cdot \alpha = (p \varepsilon(\alpha)) \otimes 
(\rho(\alpha) v), \ \ \ \alpha \in \pi.$$
Let the chain complex of $(X,\varepsilon,\rho)$ be defined as the complex of 
$\Ff[t^{\pm 1}]$-modules:
$$C^{\varepsilon,\rho}_*(X,\Ff[t^{\pm 1}]) = (\Ff[t^{\pm 1}] \otimes V) 
\otimes_{\Ff[\pi]} C_*(\widetilde{X}; \Ff).$$
It is a free based complex, where a basis is given 
by the elements of the form $1 \otimes e_i \otimes c_k$, where $\{ e_i \}$ 
is a basis of $V$ and $\{ c_k \}$ is a basis of the $\Ff[\pi]$-module
$C_*(\widetilde{X};\Ff)$, obtained by lifting cells of $X$.

A geometrical interpretation of $C_*^{\varepsilon,\rho}(X;\Ff[t^{\pm 1}])$ 
was given in~\cite{KL1}. We briefly recall their point of view.
Consider $X^\infty$ the infinite cyclic covering induced by $\varepsilon$. 
For $\overline{\pi}= \hbox{Ker } \varepsilon=\pi_1(X^\infty)$, 
the representation $\rho$ restricts to
$$ \overline{\rho} : \overline{\pi} \longrightarrow \GL(V).$$
The chain complex
$$C^{\overline{\rho}}_* (X^\infty; V) = 
V \otimes_{\Ff[\overline{\pi}]} C_*(\widetilde{X})$$
can be considered as a complex of $\Ff[t^{\pm 1}]$-modules, $\Ff[t^{\pm 1}]$ is
a trivial $\Ff[\pi]$-module, as follows:
$$t^n \cdot (v \otimes c ) = v \gamma^{-n} \otimes \gamma^n c$$
where $\gamma \in \pi$ verifies $\varepsilon(\gamma)=t$. 
In~\cite[Theorem 2.1]{KL1} it is shown that $C^{\overline{\rho}}_* (X^\infty; V)$
and $ C^{\varepsilon, \rho}_*(X,\Ff[t^{\pm 1}])$ are isomorphic as 
$\Ff[t^{\pm 1}]$-modules.

The following definition fixes some vocabulary. Denote by $\Ff(t)$ the fraction 
field of $\Ff[t^{\pm 1}]$, and define 
$C_*^{\varepsilon, \rho}(X, \Ff(t))= 
C_*^{\varepsilon,\rho}(X;\Ff[t^{\pm 1}]) \otimes_{\Ff[t^{\pm 1}]} \Ff(t)$.
 
\begin{definition}
$(X,\varepsilon,\rho)$ is \emph{acyclic} if the
chain complex $C_*^{\varepsilon, \rho}(X, \Ff(t))$ is acyclic over $\Ff(t)$.
The \emph{classical} case corresponds to the case of a \emph{trivial} $\rho$, ie.
if $V=\Ff=\Qq$ and $\rho(x)=1$ for all $x \in \pi$. 
\end{definition}

\subsection{Torsion and Alexander polynomials}

Suppose that we are given $(X,\varepsilon,\rho)$, as in the previous section.
Recall that the chain complex  $C^{\varepsilon,\rho}_*(X, \Ff(t))$ is based
by construction.

\begin{definition}
\label{torsion}
Fix a basis for the homology  $H_*^{\varepsilon, \rho}(X;\Ff(t))$.
Let $\tau_{\varepsilon, \rho}(X)$ be the torsion
of $(X, \varepsilon, \rho)$ with respect to this basis :
$$\tau_{\varepsilon, \rho}(X)=\tau(C^{\varepsilon,\rho}_*(X, \Ff(t))) \in \Ff(t)^*,$$
Up to multiplication by a factor $u t^n$ with $u \in \Ff^*$ and
$n \in \Zz$, the torsion $\tau_{\varepsilon, \rho}(X)$ is independent of the choice  
of the bases, and it is a well-defined invariant of $(X, \varepsilon, \rho)$.
\end{definition}

As mentioned in~\cite{KL1}, the indeterminacy 
of $\tau_{\varepsilon, \rho}(X)$ could be reduced. The reason for not doing so can be 
found in Theorem~\ref{homo} below.

\begin{definition}
\label{th}
The homology of $(X,\varepsilon,\rho)$ is defined as the $\Ff[t^{\pm 1}]$-module 
$$\Ho_*^{\varepsilon,\rho}(X; \Ff[t^{\pm 1}]) = 
\Ho_*(C^{\varepsilon,\rho}(X,\Ff[t^{\pm 1}])).$$
\end{definition}

One extends the definition to 
$\Ho_*^{\varepsilon,\rho}(X; \Ff(t)) = \Ho_*(C^{\varepsilon,\rho}(X,\Ff(t))$.
Since $\Ff[t^{\pm 1}]$ is a principal ideal domain, and $\Ff(t)$ is flat over 
$\Ff[t^{\pm 1}]$, one has
$$ \Ho_*^{\varepsilon,\rho}(X; \Ff(t)) \simeq 
\Ho_*^{\varepsilon,\rho}(X;\Ff[t^{\pm 1}]) \otimes \Ff(t). \eqno{(*)}$$
In particular, $\Ho_*^{\varepsilon,\rho}(X; \Ff(t))$ are $\Ff(t)$-vector spaces.

Since $\Ff[t^{\pm 1}]$ is a principal ideal domain, any $\Ff[t^{\pm 1}]$-module $H$ 
can be decomposed as a direct sum of cyclic modules. The \emph{order} of $H$ is
the product of the generators of the torsion part. If the module is free, its 
order is $1$, by convention. Note that the order of $H$ is defined up to 
multiplication $ut^n$ for $u \in \Ff^*$.

\begin{rem}
\label{rem-irred}
Note that, if $\rho=\rho_1\oplus \rho_2$ is a non-irreducible representation, then 
$$\Ho_*^{\varepsilon,\rho}(X; \Ff[t^{\pm 1}])=
\Ho_*^{\varepsilon,\rho_1}(X; \Ff[t^{\pm 1}])\oplus
\Ho_*^{\varepsilon,\rho_2}(X; \Ff[t^{\pm 1}]),$$
and hence
$$\tau_{\eps,\rho}(X)=\tau_{\eps,\rho_1}(X)\tau_{\eps,\rho_2}(X).$$
Therefore in the sequel we will only consider irreducible representations unless
otherwise stated.
\end{rem}

\begin{definition}
The $i$-th Alexander polynomial $\Delta_{\varepsilon,\rho}^i(X)$
of $(X,\varepsilon,\rho)$ is the order of 
$\Ho_i^{\varepsilon,\rho}(X;\Ff[t^{\pm 1}])$.
For short, we denote $\Delta_{\varepsilon,\rho}(X)= \frac{\Delta_{\varepsilon,\rho}^1(X)}{\Delta_{\varepsilon,\rho}^0(X)}$
the element of $\Ff(t)$.
\end{definition}

Note that in general, $\Delta_{\varepsilon,\rho}(X)$ may not be a polynomial.

\begin{thm}[\cite{Ki,KL1,Tu1}]
\label{homo}
Let $\tau_{\varepsilon, \rho}(X)$ be the torsion
of $(X, \varepsilon, \rho)$ with respect to some basis in homology. One has
$$ \tau_{\varepsilon,\rho}(X) = \frac{\prod_i \ 
\Delta_{\varepsilon,\rho}^{2i+1}(X)}{\prod_i \Delta_{\varepsilon,\rho}^{2i}(X)}.$$
\end{thm}

This theorem reconciles the two perspectives on torsion.
The first one depends on cell structure and homology bases (and as indicated above, it 
has a smaller indeterminacy). The other one depends only on the homology modules (and in particular does not require any choice of bases), and it is defined up to multiplication 
by a unit in $\Ff[t^{\pm 1}]$. For more details, see~\cite{KL1}.
Note that as a consequence, $\tau_{\varepsilon,\rho}(X)$ is a homotopy invariant.
The following lemma will be useful later.

\begin{lemma}[\cite{KL1}]
For all $(X, \varepsilon, \rho)$ such that $\varepsilon$ is non-trivial,  
$\Ho_0^{\varepsilon,\rho}(X; \Ff[t^{\pm 1}])$ is a torsion $\Ff[t^{\pm 1}]$-module.
\label{tors}
\end{lemma}

\subsection{Fox calculus} \label{Fox}
The first general definition of twisted Alexander polynomials is
due to Wada~\cite{Wa}, for finitely presented groups endowed with
an abelianization and a linear representation. Its construction involves
only Fox calculus. 

Suppose that $\pi$ has a presentation 
$\pi = \langle\ x_1, \dots ,x_m \mid r_1, \dots, r_n \ \rangle$.
The representation $\varepsilon \otimes \rho$ induces a ring 
homomorphism 
$$\array{ccc}
\Zz[\pi] & \longrightarrow &\mathcal{M}_r(\Ff[t^{\pm 1}])\\
\gamma & \longmapsto & t^{\varepsilon(\gamma)} \rho(\gamma).
\endarray$$
Let $F_m$ be the free group generated by $x_1,\dots,x_m$. Set 
$$\Phi : \Zz[F_m] \longrightarrow \Zz[\pi] 
\stackrel{\varepsilon \otimes \rho}
\longrightarrow \mathcal{M}_r(\Ff[t^{\pm 1}]).$$
Following~\cite[Lemma 1]{Wa}, there exists some $i$ such that $\Phi(x_i -1)$ 
has a non-zero determinant. Let 
$p_i : (\lambda^r)^m \longrightarrow (\lambda^r)^{m-1}$
be the projection in the direction of the $i$-th copy of  $\lambda^r$.
Consider the $(nr \times mr)$-matrix 
$$\Upsilon = \Big[ \Phi 
\big(\frac{\partial r_k}{\partial x_l} \big) \Big],$$
and define 
$$Q_i = \left\{\begin{array}{ll} \gcd \{ \big(r(m-1) \times r(m-1)\big)  
\hbox{-minors of }(p_i  \Upsilon) \} & \hbox{if } n \geq m \\ 1 & \hbox{otherwise} \end{array}\right. $$
Wada defines the twisted Alexander polynomial of $(\pi;\varepsilon,\rho)$ 
as $Q_i / \det(\Phi(x_i -1))$.
In fact, one has the following result.
\begin{thm}[\cite{KL1,Wa}] 
Let $X$ be a finite CW-complex.
If $\Ho_1^{\varepsilon,\rho}(X; \Ff[t^{\pm 1}])$ is torsion, then
$$ \Delta_{\varepsilon,\rho}(X)= \frac{Q_i}{\det(\Phi(x_i -1))}.$$
\label{foxy}
\end{thm}

\section{Duality theorem and intersection forms} \label{dua}
In this section, we recall the duality theorem for torsion, in our context.
It is due to Franz~\cite{Fr} and Milnor~\cite{Mi1}, and it
can be found in~\cite{KL1} for twisted torsion by unitary representations.

$X$ is now a compact smooth $4$-manifold, with boundary $\partial X$, possibly
empty. By Whitehead's theorem, $X$ has a canonical pl-structure, unique
up to ambient isotopy. In fact, any two pl-triangulations have a common
linear subdivision which is pl. We endow $X$ with the CW-decomposition induced
by one of these. Since $X$ is compact, the CW-complex is finite.

\begin{definition}
$(X,\varepsilon,\rho)$ is \emph{unitary} if $\rho: \pi_1(X) \longrightarrow \GL(V)$
is unitary, where $V$ is a finite dimensional $\Ff$-vector space, where $\Ff$ is a 
subfield of $\Cc$ closed under conjugation.
\end{definition}

We denote by $\overline{\cdot} : \Ff[t^{\pm 1}] \longrightarrow \Ff[t^{\pm 1}]$ 
the involution induced by the complex conjugation, and $<,>$ the positive hermitian 
form on $V$. Note that $\overline{\cdot}$ extends to $\Ff(t)$ in a natural way, so
both will be denoted the same way.

\begin{definition}
The intersection form of a unitary $(X,{\varepsilon,\rho})$ is the sesquilinear 
form
$$\array{cccl}
\varphi^{\varepsilon,\rho} : &\Ho_2^{\varepsilon,\rho} (X,  \Ff[t^{\pm 1}]) 
\times \Ho_2^{\varepsilon,\rho} (X,  \Ff[t^{\pm 1}]) &\longrightarrow&
\Ff[t^{\pm 1}]\\
&(f \otimes v \otimes c, g \otimes w \otimes c^\prime) &
\longmapsto &\sum_{\alpha \in \pi} 
( c   \cdot \alpha c^\prime) f \overline{g} \varepsilon(\alpha) 
<v \alpha, w>,
\endarray$$
where $(\cdot)$ denotes the algebraic intersection number.
\label{inter}
\end{definition}

Since for each $c$ and $c^\prime$, all but a finite number of terms are 
zero, $\varphi^{\varepsilon,\rho}$ takes values in $\Ff[t^{\pm 1}]$.

Let us fix the following convention for manifolds with boundary. We say that $(Y,\widetilde{\varepsilon},\widetilde{\rho})$ is the \emph{boundary} of $(X,\varepsilon,\rho)$ if $\partial X=Y$ as manifolds and the two following 
diagrams commute (with $\varepsilon$ and $\rho$ respectively), where
$i_*$ is induced by the inclusion:
$$\xymatrix{ \pi_1(Y) \ar[rr]^{i_*} \ar[rd]_{\varepsilon,\rho} && 
\pi_1(X) \ar[dl]^{\widetilde{\varepsilon},\widetilde{\rho} } \ar[ld] \\ & \Zz,\GL(V)}$$
From now on, $\widetilde{\varepsilon}$ and $\widetilde{\rho}$ are simply denoted
by $\varepsilon$ and $\rho$.

\begin{thm}
For any unitary $(X^4,\varepsilon,\rho)$ such that
$X$ has the homotopy type of a $2$-dimensional complex and
$C_*^{\varepsilon,\rho}(\partial X; \Ff(t))$ is acyclic, the following holds:
$$ \tau_{\varepsilon,\rho}(\partial X) = \tau_{\varepsilon,\rho}(X) \cdot
\overline{\tau}_{\varepsilon,\rho}(X) \cdot \det(\varphi^{\varepsilon,\rho}).$$
\label{dual}
\end{thm}

As mentioned before, this is a specialization of the duality theorem.
The arguments used in the proof are standard, and can be
found in particular in~\cite{KL1}. For the reader's convenience, we briefly recall it.

\begin{proof}
Note that $\partial X$ inherits the structure of a pl-manifold from $X$ and the
triangulations of $X$ can be used to define 
$C_*^{\varepsilon,\rho}(X,\partial X;\Ff(t))$, which is generated by cells
(instead of simplices). Along the lines of this proof, we will assume that the 
cell complexes have coefficients in $\Ff(t)$. Consider the pairings:
$$ C_i^{\varepsilon,\rho}(X) \times C_{4-i}^{\varepsilon,\rho}(X,\partial X)
\longrightarrow \Ff(t), \quad i=1,\dots,4,$$
which induce $\Ff(t)$-isomorphisms
$$C_{4-i}^{\varepsilon,\rho}(X,\partial X) 
\longrightarrow \overline{\hbox{Hom}}_{\Ff(t)}
(C_{i}^{\varepsilon,\rho}(X); \Ff(t)).$$
These isomorphisms take the differential of
$C_*^{\varepsilon,\rho}(X,\partial X)$ to the
dual of the differential of $ C_*^{\varepsilon,\rho}(X)$ and induce Poincare duality
isomorphisms
$$\Ho_{4-i}^{\varepsilon,\rho}(X, \partial X) 
\stackrel{\simeq} \longrightarrow H^i_{\varepsilon,\rho}(X).$$
The universal coefficient theorem applied to the chain complex 
$C_*^{\varepsilon,\rho}(X,\partial X)$
implies that evaluation induces an isomorphism
$$ H^{4-i}_{\varepsilon,\rho}(X) \simeq 
\overline{\hbox{Hom}}_{\Ff(t)} 
(\Ho_{i}^{\varepsilon,\rho}(X, \partial X) , \Ff(t)),$$
and the inner product above induces the non-singular pairing
$$\Ho_i^{\varepsilon,\rho}(X) \times 
\Ho_{4-i}^{\varepsilon,\rho}(X, \partial X) 
\longrightarrow \Ff(t). $$
For fixed bases of $ \Ho_{4-i}^{\varepsilon,\rho}(X)$, choose the dual
bases for $\Ho_i^{\varepsilon,\rho}(X, \partial X)$. We get 
$\tau _{\varepsilon,\rho}(X,\partial X)
\cdot \overline{\tau}_{\varepsilon,\rho}(X) = 1$, and by Lemma~\ref{mult}
$$ \tau_{\varepsilon,\rho}(X,\partial X) \tau_{\varepsilon,\rho}(\partial X)
\tau(\mathcal H) = \tau_{\varepsilon,\rho}(X),$$ 
where $\tau(\mathcal H)$ is the torsion of the long exact sequence in homology
of the pair $(X,\partial X)$.
Since $C_*^{\varepsilon,\rho}(\partial X; \Ff(t))$ is acyclic, $\tau(\mathcal H)$
is the alternating product of the determinant of maps induced by inclusion.
Moreover, since $X$ has the homotopy type of a $2$-complex, we have, from duality above
$$\Ho_3^{\varepsilon,\rho}(X) = \Ho_1^{\varepsilon,\rho}(X,\partial X)=0.$$
This implies in particular that the maps
$\Ho_1^{\varepsilon,\rho}(X) \rightarrow \Ho_1^{\varepsilon,\rho}(X,\partial X)$
and $\Ho_3^{\varepsilon,\rho}(X) \rightarrow \Ho_3^{\varepsilon,\rho}(X,\partial X)$
are zero. To conclude the proof, consider the diagram:
$$ \xymatrix{ \Ho_2^{\varepsilon,\rho}(X) \times 
\Ho_2^{\varepsilon,\rho}(X) \ar[rr]^{\hbox{Id} \times i_*} 
\ar[rd]_{\varphi^{\varepsilon,\rho}\otimes \Ff(t)} && 
\Ho_2^{\varepsilon,\rho}(X) \times \Ho_2^{\varepsilon,\rho}(X,\partial X) 
\ar[dl]  \ar[ld] \\ & \Ff(t) }$$

The diagonal map on the right is the unimodular map considered above. 
It follows that the matrices for the inclusion  
$\Ho_2^{\varepsilon,\rho}(X) \rightarrow \Ho_2^{\varepsilon,\rho}(X,\partial X)$ 
and $\varphi^{\varepsilon,\rho}\otimes \Ff(t)$ are conjugated. Hence, 
$\tau(\mathcal H) = 1 / \det (\varphi^{\varepsilon,\rho}).$
Note that since $X$ has the homotopy type of a $2$-complex, 
$H_2 ^{\varepsilon,\rho}(X;\Ff[t^{\pm 1}])$
is free and a matrix for $\varphi^{\varepsilon,\rho}$ is also a matrix for 
$\varphi^{\varepsilon,\rho} \otimes \Ff(t)$.
\end{proof}

\section{Knots and Links}
\label{hist}
In this section, we briefly collect results on the torsion applied to 
link complements and Alexander polynomials.
Let $L$ be an oriented link in $S^3$ with $\mu$ components,
and $E$ be the exterior of $L$. Denote the homology classes of the meridians of the
link components by $\nu_i$. Note that
$$\Ho_1(E) = \oplus_{i=1}^\mu \Zz \nu_i.$$
Let $q_1,\dots,q_\mu$ be integers with $\gcd(q_1,\dots,q_\mu)=1$, and define
$$\varepsilon : \Ho_1(E) \longrightarrow \Zz= \langle\ t\ \rangle$$
$$\nu_i \longmapsto t^{q_i}.$$
Since $\gcd(q_1,\dots,q_\mu)=1$, the associated infinite cyclic covering is connected.
Let $\rho: \pi_1(E) \longrightarrow \GL(V)$ for some finite dimensional $V$ over a 
field~$\Ff$.

\begin{definition}
The twisted torsion of $(L,\varepsilon,\rho)$ is 
$\tau_{\varepsilon,\rho}(L) = \tau_{\varepsilon,\rho} (E)$.
Similarly, the twisted Alexander polynomials are 
$\Delta_{\varepsilon,\rho}^i(L)= \Delta_{\varepsilon,\rho}^i(E)$, and denote $\Delta_{\varepsilon,\rho}(L)= \Delta_{\varepsilon,\rho}^1(L)
/ \Delta_{\varepsilon,\rho}^0(L)$.
\label{deflink}
\end{definition}

For knots, Cha~\cite{Cha} also considered the case where $\Ff$ is not a field 
but any Noetherian factorization domain. If $\rho$ is the trivial representation 
and $q_i=1$ for all $i$, $\Delta_{\varepsilon,\rho}^1(L)=\Delta_L$ is the classical
Alexander polynomial.

A CW-structure on $E$ can be given with one $0$-cell, $n$ $1$-cells and $(n-1)$
$2$-cells. This can correspond to a Wirtinger presentation 
$\pi_1(E)= \langle\ x_1,\dots,x_n \mid r_1,\dots, r_{n-1}\ \rangle$,
where $r_i$ is of the form $x_j x_k x_j^{-1} x_l^{-1}$. From this, we deduce that
the chain complex of $(E,\varepsilon,\rho)$ is given by the following sequence 
$$0 \longrightarrow (\Ff[t^{\pm 1}] \otimes_\Ff V)^{n-1} 
\stackrel{\partial_2} \longrightarrow
(\Ff[t^{\pm 1}] \otimes_\Ff V)^{n} \stackrel{\partial_1} \longrightarrow \Ff[t^{\pm 1}]
\otimes_\Ff V \longrightarrow 0.$$
A matrix for $\partial_2$ is given by the image of Fox derivatives, see Section 2.4. 
The map $\partial_1$ has a column matrix with entries $\rho(x_i)t^{\eps(x_i)} - \Id$.
Note that this complex may not be acyclic in general, even for knots (see for 
example~\cite[Section 3.3]{KL1}). In any case, one has 
$\Delta_{\varepsilon,\rho}^2(L)=1$ and hence the following proposition is a direct
consequence of Theorem~\ref{homo}.

\begin{prop}
For any oriented link $L$ in $S^3$, and $(\varepsilon,\rho)$, one has
$$\tau_{\varepsilon,\rho}(L) = \Delta_{\varepsilon,\rho}(L).$$
In particular, if $\rho$ is the trivial representation and $q_i=1$ for all $i$,
then 
$$(t-1)\tau_{\varepsilon,\rho}(L) = \Delta_L.$$
\label{alexlink}
\end{prop}

Note that it is easy to compute $\Delta^0_{\varepsilon,\rho}(L)$ from
 the description of $\partial_1$ above.

Proposition~\ref{alexlink} was first shown by Milnor~\cite{Mi2} 
(with several variables) and extended by Turaev to any $3$-manifold.
In the twisted case, this is due to Kitano for knot complements~\cite{Ki}
and later extended to links by Wada~\cite{Wa}, in the acyclic case.
It is worth mentioning that Cha~\cite{Cha} gave a formula for  
$\Delta_{\varepsilon,\rho}(K)$ when $K$ is a fibered knot, in terms of the 
homotopy type of the monodromy.

\section{Plane curves}
\label{planecurves}
Let $\mathcal C$ be an algebraic curve in $\Cc^2$ with $r$ 
irreducible components. Let us denote by $X$ the complement in $\Bb^4$ of an open
tubular neighborhood of $\mathcal C$, for a sufficiently large
ball $\Bb^4 \subset \Cc^2$. As in the case of links 
we say $\gamma_\ell\in \pi_1(X)$ is a meridian of $\cC_\ell$ if 
$\gamma_\ell=\beta \cdot \tilde{\gamma}_\ell \cdot \beta^{-1}$, where $\beta$
is a path from the base point to a point on the boundary of a tubular neighborhood 
$T_\ell$ of $\cC_\ell$ and $\tilde{\gamma}_\ell$ is a positively oriented loop that 
bounds a fiber of $T_\ell$. It is well known that $\pi_1(X)$ is generated by 
meridians of the irreducible components. Hence $\Ho_1(X) = \Zz^r$ is generated 
by the homology classes $\nu_\ell$ of the meridians $\gamma_\ell$ of $\mathcal C_\ell$
for $\ell=1,\dots,r$.

\subsection{Twisted torsion of plane curves} \label{plc}
Let $q_1,\dots,q_r \neq 0$ be integers with $\gcd(q_1,\dots,q_r)=1$.
Consider
$$\varepsilon : \Ho_1(X) \longrightarrow \Zz$$
with $\epsilon(\nu_\ell)=q_\ell$. Let
$$\rho: \pi_1(X) \longrightarrow \GL(V),$$
be a fixed representation where $V$ is a finite dimensional $\Ff$-vector space.

\begin{definition}
The torsion of $(\mathcal C,\varepsilon,\rho)$ is 
$\tau_{\varepsilon,\rho}(\mathcal C)=\tau_{\varepsilon,\rho}(X)$. Similarly,
Alexander polynomials can be defined as
$\Delta^i_{\varepsilon,\rho}(\mathcal C)=\Delta^i_{\varepsilon,\rho}(X)$.
As in the case of links, we will denote
$\Delta^1_{\varepsilon,\rho}(\mathcal C) / \Delta^0_{\varepsilon,\rho}(\mathcal C)$ 
simply by $\Delta_{\varepsilon,\rho}(\mathcal C)$.
\end{definition}

In the $1930$'s Zariski and Van-Kampen developed a method to compute the fundamental 
group of $X$. Refinements of this algorithm were constructed later mainly by Chisini 
and Moishezon. Finally, Libgober described the homotopy type of $X$ as follows. 

Consider a generic linear projection $\pi:\Cc^2 \rightarrow \Cc$, i.e. such that: 
\begin{enumerate}
\item there are no vertical asymptotes,
\item the fibers are transversal to $\mathcal C$ except for a finite number of them 
which are either simple tangents to a point of $\cC$ or lines through a singular point 
of $\cC$ transversal to its tangent cone.
\end{enumerate}
Let $\mathcal P$ be the (finite) set of critical values of $\pi$. 
The braid monodromy of $\cC$ is the homomorphism
$$\vartheta : \pi_1(\Cc\setminus \mathcal P) \longrightarrow \Bb_d,$$ 
where $\Bb_d$ denotes the braid group, viewed as the mapping class group of a generic
fiber relative to $\cC$, that is, $(\pi^{-1}(p),\pi^{-1}(p)\cap \cC)$ with 
$p\in \Cc\setminus \mathcal P$. We fix a basis $\{ \alpha_i \}_{i=1,...,n}$ 
of $\pi_1(\Cc \setminus \mathcal P)$ such that:

\begin{enumerate}
\enet{(B\arabic{enumi})}
\item
\label{bmc1}
$\alpha_i=A_i\cdot s_i \cdot A^{-1}$ where
$\bar A_i$ is a path joining $p$ and $p_i\in \mathcal P$, $s_i$ is the boundary
of a small disc $D_i$ around $p_i$ and $A_i=\bar A_i \setminus D_i$ 
\item
\label{bmc2}
$\alpha_1\cdot \ldots \cdot \alpha_n$ is 
homotopic to the boundary of a big disc containing $\mathcal P$. 
\end{enumerate}
Also fix a basis $\gamma_1,\dots,\gamma_d$ of the (free) fundamental group 
$\pi_1(\pi^{-1}(p)\setminus \cC)$. 

Note that each $\pi^{-1}(\bar A_i)$ produces the collapse of $m_i$ points of
$\pi^{-1}(p)\cap \cC=\{q_1,...,q_d\}$ to a point on $\cC$, say $P_i$, where $m_i$ 
is nothing but the multiplicity of intersection of the line $\pi^{-1}(p)$ and 
the curve $\cC$ at $P_i$. We will denote such points by 
$q_{i,1},...,q_{i,m_i}\in \pi^{-1}(p)\cap \cC$ and analogously their meridians
$\gamma_{i,1},...,\gamma_{i,m_i}$. The action of $\vartheta(\alpha_i)$ on 
$\gamma_{i,k_i}$, $k_i=1,...,m_i$ can be decomposed via $A_i$ and $s_i$ as follows 
$$\vartheta(\alpha_i)(\gamma_{i,k_i})=
\sigma_i(\omega_{i,k_i} \gamma_{i,k_i} \omega_{i,k_i}^{-1})$$
where $\sigma_i$ only depends on the local type of the singularity of the projection
and $\tilde \gamma_{i,k_i}=\omega_{i,k_i} \gamma_{i,k_i} \omega_{i,k_i}^{-1}$ is 
homotopic to a meridian of $\cC$ on $S^3_i \setminus \cC$ (a small 3-dimensional 
sphere centered at $P_i$).

In the following theorem, the presentation of $\pi_1(X)$ is due to~\cite{Za,VK}. 
Libgober~\cite{Li2} used this presentation to describe the homotopy type of $X$.
We recall this well-known result in order to set notation for future reference.

\begin{thm} 
\label{homtype}
The two-dimensional complex associated with the following presentation of $\pi_1(X)$:
\begin{equation}
\label{eq-pres}
\langle\ \gamma_1,\dots,\gamma_d \mid 
\sigma_i(\tilde \gamma_{i,k_i})=\tilde \gamma_{i,k_i}, \ 
i=1,\dots,n, \ k_i=1,\dots,m_i-1\ \rangle,
\end{equation}
has the homotopy type of $X$.
\end{thm}

Note that Theorem~\ref{homtype} provides in particular a presentation of the fundamental
group and can be useful to construct $\rho$ explicitly. 

\begin{prop}
For any algebraic curve $\mathcal C$ in $\Cc^2$, one has
$$ \tau_{\varepsilon,\rho}(\mathcal C) = \Delta_{\varepsilon,\rho}(\mathcal C).$$
In particular, if $\rho$ is the trivial representation and $q_i=1$, then
$$(t-1)\tau_{\varepsilon,\rho}(\mathcal C) = \Delta_\mathcal C.$$
\label{alexcurve}
\end{prop}

In fact, the second statement holds for any epimorphism $\eps$ and 
$\tau_{\varepsilon}(\mathcal C)$ coincides with the corresponding $\theta$-polynomial introduced by Oka~\cite{Oka2}.
Note that the complex $C_*^{\varepsilon,\rho}(X; \Ff(t))$ is not acyclic in general. 
In fact even the Euler characteristic is not zero.

Similarly to the case of links 
$\Delta^0_{\eps,\rho}(\cC)$
 can be easily computed.

\begin{proof}
By Theorem~\ref{homtype}, $X$ has the homotopy type of a $2$-complex.
From this --\,and similarly to link complements\,--, we deduce that
the chain complex of $(X,\varepsilon,\rho)$ is given by the following sequence: 
$$ 0 \longrightarrow (\Ff[t^{\pm 1}] \otimes_\Ff V)^{n-1} \stackrel{\partial_2} 
\longrightarrow (\Ff[t^{\pm 1}] \otimes_\Ff V)^{n} \stackrel{\partial_1} \longrightarrow
\Ff[t^{\pm 1}] \otimes_\Ff V \longrightarrow 0.$$
Clearly $H_i^{\varepsilon,\rho}(X; \Ff[t^{\pm 1}])=0$ for $i \geq 3$.
Since,  $H_2^{\varepsilon,\rho}(X; \Ff[t^{\pm 1}])$ is free or zero,
$ \Delta^2_{\varepsilon,\rho}(\mathcal C) =1$.
Hence by Theorem~\ref{homo},
$$ \tau_{\varepsilon,\rho}(\mathcal C) = \frac{\Delta^1_{\varepsilon,\rho}
(\mathcal C)}{\Delta^0_{\varepsilon,\rho}(\mathcal C) \cdot \Delta^2_{\varepsilon,\rho}(\mathcal C)}.$$
Moreover, if $\rho$ is trivial, 
$$\Delta^0_{\varepsilon}(\mathcal C)=\gcd (t^{q_\ell} - 1)=t-1.\qedhere$$
\end{proof}

In some cases, one can assure that $\Delta_{\eps,\rho}(\cC)$ is actually a polynomial.

\begin{prop}
If $\cC$ is not irreducible and $H^{\varepsilon,\rho}_1(X;\Ff[t^{\pm 1}])$ is torsion, 
then $\Delta_{\eps,\rho}(\cC)$ is a polynomial. 
\end{prop}

\begin{proof}
Consider $\text{ab}:H_1(\Cc^2\setminus \cC) \to \Zz^r$ the abelianization morphism.
We will consider a several-variable invariant defined by Wada~\cite{Wa} as 
$\frac{Q_i(t_1,\dots,t_r)}{P_i(t_1,\dots,t_r)}$, where $Q_i$ is the $\gcd$ of the minors
of maximal order of $p_i\Upsilon_{\text{ab},\rho}$ and $P_i=\det(\Id-\rho(\nu_i)t_i)$.
Following the proof of \cite[Proposition~9]{Wa},
$\frac{Q_i(t_1,\dots,t_r)}{P_i(t_1,\dots,t_r)}$ is a
polynomial in the variables $(t_1,\dots,t_r)$. 
Since $H^{\varepsilon,\rho}_1(X;\Ff[t^{\pm 1}])$ is torsion, one has
$\Delta_{\eps,\rho}(\cC)=\frac{\tilde Q_i(t)}{\tilde P_i(t)}$ (Theorem~\ref{foxy}),
where $\tilde Q_i$ is the $\gcd$ of the minors of maximal order of 
$p_i\Upsilon_{\eps,\rho}$ and $\tilde P_i=\det(\Id-\rho(\nu_i)t^{q_i})$. Note that
$P_i(t^{q_1},...,t^{q_r})=\tilde P_i(t)$ and $P_i(t_1,...,t_r)$ divides every minor in 
$p_i\Upsilon_{\text{ab},\rho}$. Therefore $P_i(t^{q_1},...,t^{q_r})=\tilde P_i(t)$ 
divides every minor of $p_i\Upsilon_{\eps,\rho}$. This means that $\tilde P_i(t)$ also
divides $\tilde Q_i(t)$, and thus $\Delta_{\eps,\rho}(\cC)$ is also a polynomial.
\end{proof}

\subsection{Relation with local polynomials}

Suppose that we are given $(\mathcal C,\varepsilon, \rho)$.
Let $S_1^3,\dots,S_s^3$ be sufficiently small $3$-spheres around the singular points 
$\{P_1,...,P_s\}$ of $\cC$. Denote by $L_k= \mathcal C \cap S^3_k$ the link of the
singularity at $P_k$, and by $E_k$ be the link exterior. Also choose a base point
$Q_i\in S^3_k \setminus L_k$ and denote by 
$\pi^k=\pi_1(S^3_k \setminus L_k;Q_i)$ the local fundamental groups at $P_k$. 
The inclusion maps 
$i_k : \pi^k \rightarrow \pi_1(X)$ and $(\varepsilon,\rho)$ induce morphisms
$$\varepsilon_k: \pi^k \longrightarrow \Zz $$
$$\hbox{and } \rho_k: \pi^k \longrightarrow \GL(V),$$
for any $k=1,...,s$.

\begin{definition}
For all $k=1,\dots,s$, let the local torsions be defined as
$$\tau_k = \tau_{\varepsilon_k,\rho_k}(L_k).$$
Analogously, by Theorem~\ref{alexlink}, we define
$\Delta_k=\Delta_{\varepsilon_k,\rho_k}(L_k)=\tau_k$.
Also, there is a local torsion \emph{at infinity} defined as
$$\tau_\infty = \tau_{\varepsilon_\infty,\rho_\infty}(L_\infty),$$
where $L_\infty=\cC \cap \partial \Bb^4$ is the intersection of the
curve with the boundary of the big ball considered at the beginning
of the section.

\end{definition}

Note that an explicit description of the maps $\pi^k \rightarrow \pi$ can be 
obtained via the braid monodromy of a generic projection of the curve (accurate
packages have been developed to compute braid monodromies of curves with
equations over the rationals by Bessis~\cite{Bes} and Carmona~\cite{Ca}). 
For notation the reader is referred to the discussion previous to 
Theorem~\ref{homtype}(\ref{eq-pres}). Note that, via the lifting of the paths $A_k$, 
the meridians $\tilde \gamma_{k,j}$ can also be seen as meridians of $\cC$ on 
$S^3_k$ based at $Q_k\in S^3_k \setminus L_k$. The local fundamental group $\pi^k$
can be presented as
$$\pi^k=
\langle\ \tilde \gamma_{k,j} \mid \sigma_k(\tilde \gamma_{k,j})=
\tilde \gamma_{k,j}, j=1,...,m_k-1 \ \rangle$$
and $\rho_k(\tilde \gamma_{k,j})=\rho(\omega_{k,j} \gamma_{k,j} \omega_{k,j}^{-1})$, 
for any $j=1,...,m_k-1$.
A similar description of $\pi^\infty,\varepsilon_\infty$, and $\rho_\infty$ can
be given analogously using $\vartheta(\alpha_1\cdot ...\cdot \alpha_n)$ 
(as defined in~\ref{bmc2}) instead of the local braids.

Let $\varphi^{\varepsilon,\rho}(\cC)$ be the intersection form of $(X,\varepsilon,\rho)$
given in Definition~\ref{inter}. Theorem~\ref{thm-main-intro}, stated in the introduction, 
can be written as follows.

\begin{thm}
Let $(\mathcal C,\varepsilon, \rho)$ be unitary and suppose that the local 
representations $\rho_k$ are acyclic. Then
$$\Big( \prod_{\ell=1}^r \det(\Id - \rho(\nu_\ell) 
t^{q_\ell})^{s_\ell-\chi(\mathcal C_\ell)}\Big) \cdot \prod_{k=1,...,s,\infty} 
\Delta_k = \Delta_{\varepsilon,\rho}(\mathcal C) \cdot
\overline{\Delta}_{\varepsilon,\rho}(\mathcal C) \cdot 
\det \varphi^{\varepsilon,\rho}(\cC),$$
where $\nu_\ell$ is the homology class of a meridian around the irreducible 
component $\cC_\ell$ and $s_\ell=\# \Sing(\cC) \cap \cC_\ell$.
\label{main}
\end{thm}

\begin{rem}
\label{rem-nu}
From the discusion above, since two meridians of the same irreducible component
are conjugated, $\det(\Id-\rho(\gamma_\ell)t^q)$ only depends on the homology class 
of $\gamma_\ell$, say $\nu_\ell$. By abuse of notation, we simply write
$\det(\Id-\rho(\nu_\ell)t^q)$ for this element of $\Ff[t^{\pm 1}]$.
\end{rem}

\begin{cor}
If $\varphi^t(\cC)$ is the intersection form with
twisted coefficients in $\Qq[t^{\pm 1}]$, then
$$(t-1)^{1-\chi(\cC)}\prod_{k=1,...,s,\infty} \Delta_{L_k} = \Delta_{\mathcal C}^2
\cdot \det \varphi^t(\cC).$$ 
\label{corro}
\end{cor}

\begin{proof}
An immediate consequence of Propositions~\ref{alexlink},~\ref{alexcurve}, and
Theorem~\ref{main} considering $\rho$ trivial and $q_i=1$, $i=1,...,r$; 
in addition to the obvious relation 
$\sum_{\ell=1}^r \left( s_\ell - \chi(\cC_\ell) \right)=s - \chi(\cC)$,
and to the fact that $\overline{\Delta}_{\mathcal C}=\Delta_{\mathcal C}$ 
(since $\Delta_{\mathcal C}$ is a product of cyclotomic polynomials).
\end{proof}

\begin{proof}[\textit{Proof of Theorem~\ref{main}.}]

Let $(M,\varepsilon,\rho)$ be the boundary of the curve exterior and let
$F= \mathcal C \setminus \sqcup_k (\mathcal C \cap \Bb^4_k)$, $k=1,...,s,\infty$
be the (abstract) surface obtained by removing disks $D_1 \cup \dots \cup D_{n_k}$ 
from the boundary of a ball around each singular point $P_k$ of $\mathcal C$, 
$k=1,...,s$ and from the boundary of the big ball $\Bb^4=\Bb^4_\infty$. 
Let $N= F \times S^1$. The boundary of $N$ is a union of disjoint tori 
$T_1^k \cup \dots \cup T_{n_k}^k$ for $k=1,\dots,s,\infty$. From a plumbing 
description of a tubular neighborhood of the curve, one can show that $M$ is obtained 
by gluing the link exteriors $E_k$ with $N= F \times S^1$, along the tori $T^k_i$ for
$i=1,\dots,n_k$, $k=1,...,s$: 
$$M = N \cup_{\coprod_i T_i^k} \big( \coprod S^3_k \setminus L_k \big).$$ 
The gluing map sends a longitude of $L_k^i$ to the restriction
of a section in $N$, and a meridian of $L_k^i$ to a fiber in $N$.
The inclusion maps induce triples $(N,\varepsilon,\rho)$, 
$(S^3_k \setminus L_k,\varepsilon_k,\rho_k)$, and $(T_i^k,\varepsilon_k,\rho_k)$ 
(with $k=1,\dots,s,\infty$, $i=1,\dots,n_k$).
By multiplicativity, $\tau_{\varepsilon,\rho}(N)$ is the product of the torsion of the connected components of $N$. Let us compute $(F_\ell \times S^1,\varepsilon,\rho)$ where
$F_\ell$ is a connected component of $F$, corresponding to $\cC_\ell$. Note that 
$H_1(F_\ell \times S^1)= H_1(F_\ell) \oplus \Zz \nu_\ell$ and
$\varepsilon(\nu_\ell)=q_\ell$. Hence by Lemma~\ref{tors},
$H_0^{\varepsilon,\rho}(F_\ell \times S^1)$ is torsion over $\Ff[t^{\pm 1}]$.
If $F_\ell$ is a disk, then $\pi_1(F_\ell \times S^1)=\Zz$
and $\tau_{\varepsilon,\rho}(F_\ell \times S^1)=1$. If $F_\ell$ is not a disk, then
$\pi_1(F_\ell \times S^1)$ can be presented by $b_1(F_\ell) + 1$ generators 
$a_1,\dots,a_{b_1(F_\ell)},\gamma_\ell$ (such that $\bar \gamma_\ell=\nu_\ell$) 
with the relations $a_k \gamma_\ell = \gamma_\ell a_k$ for
all $k=1,\dots,b_1(F_\ell)$. Taking the matrix of Fox derivatives of the relations and 
tensoring with $\Ff[t^{\pm 1}] \otimes V$ one obtains a 
$b_1(F_\ell) \times (b_1(F_\ell)+1)$ matrix, with entries in 
$M_{m \times m}(\Ff[t^{\pm 1}])$ ($m=\dim_{\Ff}V$) describing the differential
$$\partial_2: C_2(F_\ell \times S^1; \Ff[t^{\pm 1}]) \longrightarrow 
C_1(F_\ell \times S^1; \Ff[t^{\pm 1}]).$$
The $i$-th row of this matrix has $\Id - \rho(\gamma_\ell) t^{q_\ell}$ in the 
$\ell$-th column and $\Id - \rho(a_i) t^{\varepsilon(a_i)}$ in the last column.
The zero chains are $C_0(F_\ell \times S^1; \Ff[t^{\pm 1}])=\Ff[t^{\pm 1}] \otimes V$
and $\partial_1$ is the column matrix with $i$-th entry 
$\rho(a_i) t^{\varepsilon(a_i)} - \Id$ for $i \leq b_1(F_\ell)$ and last
entry $\rho(\gamma_\ell) t^{q_\ell}- \Id$.
Dropping the last column of the matrix for $\partial_2$,
one obtains a $b_1(F_\ell) \times b_1(F_\ell)$ matrix with determinant
$\det(\Id - \rho(\gamma_\ell) t^{q_\ell})^{b_1(F_\ell)}$. Note that this determinant
only depends on the conjugacy class $\nu_\ell$ of $\gamma_\ell$, and hence, it 
does not depend on the chosen meridian of $\cC_\ell$. Therefore, in the future we
will simply write $\det(\Id - \rho(\nu_\ell) t^{q_\ell})$. Also note that the chain complex
$C_*^{\varepsilon,\rho}(F_\ell \times S^1; \Ff(t))$ is acyclic since $q_\ell \neq 0$.
From Theorem~\ref{foxy}, one obtains
$\tau_{\varepsilon,\rho}(F_\ell \times S^1) =
\det(\Id - \rho(\nu_\ell) t^{q_\ell})^{-\chi(F_\ell)}$.
Note that $\chi(F_\ell)=\chi(\cC_\ell) - s_\ell$.
By additivity of Euler characteristic, one has
$$\tau_{\varepsilon,\rho}(N) = \prod_{\ell=1}^r
\det(\Id - \rho(\nu_\ell) t^{q_\ell})^{s_\ell-\chi(\mathcal C_\ell)}.$$
Also note that, since $\pi_1(T_i^k)$ is Abelian, one has
$$\prod_{k=1,...,s,\infty}\prod_{i=1}^{n_k} \tau_{\varepsilon_k,\rho_k}(T_i^k) =1.$$

Hence one has the following Mayer-Vietoris sequence with coefficients in $\Ff(t)$:
$$0 \longrightarrow \oplus_{i,k} C_*^{\varepsilon_k,\rho_k}(T_i^k)
\longrightarrow \oplus_k \big(C_*^{\varepsilon_k,\rho_k}(S^3_k \setminus L_k)) \oplus
C_*^{\varepsilon,\rho}(N) \longrightarrow C_*^{\varepsilon,\rho}(M) \longrightarrow 0.$$
By hypothesis the chain complex 
$\oplus_k C_*^{\varepsilon_k,\rho_k}(S^3_k\setminus L_k;\Ff(t))$ 
is acyclic. Since  $C_*^{\varepsilon,\rho}(N; \Ff(t))$ and
$\oplus_{i,k} C_*^{\varepsilon_k,\rho_k}(T_i^k)$ are
also by computations above, then $C_*^{\varepsilon,\rho}(M)$ is acyclic.
From Lemma~\ref{mult}, we obtain
$$\big( \prod_{k} \tau_k \big) \cdot \tau_{\varepsilon,\rho}(N)=
\big( \prod_{k,i} \tau_{\varepsilon_k,\rho_k}(T_i^k) \big)
\cdot \tau_{\varepsilon,\rho}(M).$$
Hence,
$$\tau_{\varepsilon,\rho}(M)= \Big( \prod_{\ell=1}^r \det(\Id - \rho(\nu_\ell) t^{q_\ell})^{s_\ell-\chi(\mathcal C_\ell)} \Big) \cdot \prod_{k} \tau_k.$$
By Theorems~\ref{dual} and~\ref{alexcurve}, one obtains the result.
\end{proof}

\subsection{Twisted Alexander Polynomials and Characteristic Varieties}
Characteristic Varieties were introduced by Hillman~\cite{Hi} for links and
by Libgober~\cite{Li4} for algebraic curves. They can be defined as the zero set 
of the first Fitting ideal $F_1(X)$ of $H_1({X_{\text{ab}}})$ as a 
$\Cc[H_1(X)]$-module (where ${X_{\text{ab}}}$ represents the universal Abelian
cover of $X$). We will briefly describe the close relationship between Characteristic 
Varieties and twisted Alexander invariants associated with rank one representations.

Let $\bar \xi=(\xi_1,...,\xi_r)$ be in the zero set of the first characteristic variety
$V_1(X)\subset (\Cc^*)^r$ associated with an affine curve $\cC\subset \Cc^2$
with $r$ irreducible components $\cC_1$, $\cC_2$,..., $\cC_r$. Consider the 
representation $\rho:\pi_1(X)\to \Cc^*=\GL(1,\Cc)$ given by
$\rho_{\bar \xi}(\gamma_i)=\xi_i$, for any $\gamma_i$ meridian of the irreducible
component $\cC_i$ and $\eps:H_1(X)\to \Zz$ such that $\eps(\nu_\ell)=1$. Hence
any polynomial $p(t_1,...,t_r)\in F_1(X)$ satisfies $p(\bar \xi)=0$. Note that
$\Delta^1_{\eps,\rho_{\bar \xi}}(\cC)=\gcd\{p(\xi_1 t,...,\xi_r t) \mid p\in F_1(X)\}$.
Hence, either $p(\xi_1 t,...,\xi_r t)=0$ for all $p\in F_1(X)$ or 
$\Delta^1_{\eps,\rho_{\bar \xi}}(\cC)$ contains $(t-1)$ as a factor, since 
$t=1$ is a root of $p(\xi_1 t,...,\xi_r t)$. Therefore, $\bar \xi\in V_1(X)$ 
if and only if $(t-1)$ divides $\Delta^1_{\eps,\rho_{\bar \xi}}(\cC)$.
Hence, one has the following result.

\begin{prop}
There is a one-to-one correspondence between the first characteristic variety of $\cC$ 
and the set of 1-dimensional representations of $\pi_1(X)$ with non-trivial twisted
Alexander polynomial $\Delta^1_{\eps,\rho}(\cC)$.
\end{prop}

Note that here \emph{non-trivial} means \emph{not a unit}, that is, zero is considered
a non-trivial polynomial.

\section{Examples}
\label{exs}

\subsection{A Zariski pair} Consider the space ${\mathcal M}$ of sextics with the
following combinatorics:
\begin{enumerate}
\smallbreak\item $\cC$ is a union of a smooth conic $\cC_2$
and a quartic~$\cC_4$. 
\smallbreak\item $\Sing(\cC_4)=\{P,Q\}$ where $Q$ is a cusp of
type $\Aa_4$ and $P$ is a node of type~$\Aa_1$.
\smallbreak\item $\cC_2 \cap \cC_4=\{Q,R\}$ where $Q$ is a 
$\Dd_7$ on $\cC$ and $R$ is a $\Aa_{11}$ on~$\cC$.
\end{enumerate}
As shown in~\cite{ACC}, ${\mathcal M}$ consists of two irreducible components.
One can add a transversal line and calculate the fundamental groups of representatives
$\cC^{(1)}$ and $\cC^{(2)}$ in each component. One has the following:
$$\pi_1(\Cc^2\setminus \cC^{(1)})=\langle e_1,e_2 : [e_2,e_1^2]=1,
(e_1e_2)^2=(e_2e_1)^2, [e_1,e_2^2]=1 \rangle$$
$$\pi_1(\Cc^2\setminus \cC^{(2)})=\langle e_1,e_2 : [e_2,e_1^2]=1,
(e_1e_2)^2=(e_2e_1)^2 \rangle.$$
Assuming $\eps:\Zz^2\to \Zz$, $\eps(1,0)=\eps(0,1)=1$ and taking, for instance, 
the rank~1 representation $\rho(e_1)=1$, $\rho(e_2)=-1$ (which is unique up to 
equivalence) one obtains:
$$\Delta_{\eps,\rho}(\cC^{(1)})=1,\quad \Delta_{\eps,\rho}(\cC^{(2)})=t+1.$$
Also, note that the classical Alexander polynomial is trivial for both curves.

\subsection{Nodal degenerations}
The following example illustrates how twisted Alexander polynomials can be 
sensitive to nodes, something that classical Alexander polynomials are not.

We say a curve $\cD$ is of type $I$ if $\cD$ is an irreducible plane curve of 
degree $d$ such that $\cD$ has an ordinary $(d-2)$-ple point at $P$.
Let $L_1$ and $L_2$ be lines through $P$ such that either $L_i$ is tangent to a smooth
point $P_i\in \cD$ or $L_i$ passes through a double point $P_i\neq P$ of type $\Aa_{2r}$.
Let us denote $\cC = L_1 + L_2 + \cD$. Assume that $\cD$ has only nodes as singular points
apart from $P$. According to~\cite[Theorem 1]{ACT}, $\cD$ is rational 
if and only if there exists a dihedral cover $\mathbb D_{2n}$ ramified along 
$2(L_1+L_2)+n\cD$ for any $n\geq 3$. In fact, according to~\cite[Corollary 2]{ACT2} 
this implies that $\cD$ is rational if and only if the fundamental group of 
$\Pp^2\setminus (L_1\cup L_2\cup \cD)$ admits $\Zz_2*\Zz_2$ as a quotient. 
Moreover (see~\cite[Proposition 6.1]{ACT}), there exist nodal degenerations 
$\cD_t\to \cD_0$ to a rational $\cD_0$ of type $I$ using (non-rational) curves 
$\cD_t$ ($t> 0$) of type $I$ with Abelian fundamental groups. A presentation for the
fundamental group of $\cC_0=L_1 + L_2 + \cD_0$ is given as follows:
$$
G(\cC_0)=\!\!\langle\ \ell,x_1,x_2\mid
[x_1,x_2]=1, 
\ell^{-1} x_1 \ell= x_2,\ell^{-1} x_2 \ell= x_1 \ \rangle.$$ 
Considering $\eps$ the usual morphism $\eps(\nu_\ell)=1$, and 
$$\rho(\ell)=\left(\array{cc} 1&0\\0&-1\endarray\right),$$
$$\rho(x_1)=\left(\array{cc} -1&0\\1&-1\endarray\right),$$
$$\rho(x_2)=\left(\array{cc} -1&0\\-1&-1\endarray\right)$$
one obtains 
\begin{equation}
\label{eq-tap}
\Delta_{\eps,\rho}(\cC_0)=t+1.
\end{equation}
Note that $\rho(G(\cC_0))\cong \Zz_2*\Zz_2$.

The three non-nodal singularities of $\cC_t$ are $\{P,P_1,P_2\}$
and they lie on the lines $L_1$ and $L_2$. Hence, maybe by performing
projective transformations, we can assume that $\{P,P_1,P_2\}$ and
$L_1$ and $L_2$ are fixed throughout the degeneration. This implies
that the classical Alexander polynomial $\Delta_{\cC_t}$ of $\cC_t$ is 
the same for all $t\geq 0$ (since they have the same adjunction ideals,
see for instance~\cite[Theorem 5.1]{Lib83}). Since $G(\cC_t)$ is Abelian, 
this implies that $\Delta_{\cC_1}=\Delta_{\cC_0}=1$. Formula~(\ref{eq-tap}) 
shows that $\cC_0$ has a non-trivial twisted Alexander polynomial.

\end{document}